\documentclass[12pt,reqno]{amsart}

\usepackage[utf8]{inputenc}

\usepackage{amsthm}

\theoremstyle{plain}
\newtheorem{theorem}{Theorem}
\newtheorem{remark}{Remark}
\newtheorem{proposition}{Proposition}[section]

\newtheorem{corollary}[theorem]{Corollary}

\numberwithin{equation}{section}

\begin{document}

\title[Beltrami equation for the harmonic diffeomorphisms]{Beltrami equation for the harmonic diffeomorphisms between surfaces
}

\author{A. Fotiadis}

\email{A. Fotiadis:  fotiadisanestis@math.auth.gr}
\author{C. Daskaloyannis}
\email{C. Daskaloyannis: daskalo@math.auth.gr }
\curraddr{
	Department of Mathematics, Aristotle University of Thessaloniki,
Thessaloniki 54124, Greece}
\date{}
\subjclass[2010]{58E20, 53C43, 35Q53, 30C62}

\keywords{harmonic maps, Beltrami equations, sinh-Gordon and sin-Gordon equations. }

\begin{abstract}
In references \cite{Minsky,Wolf1} it was proved that harmonic diffeomorphisms, with nonvanishing Hopf differential, satisfy a Beltrami equation of a certain type: the imaginary part of the logarithm of the Beltrami coefficient coincides with the imaginary part of the logarithm of the Hopf differential. Therefore, it is a harmonic function. The real part of the logarithm of the Beltrami coefficient satisfies an elliptic nonlinear differential equation, which in the case of constant curvature is the elliptic sinh-Gordon equation.  

In this paper we also prove the converse: if the imaginary part of the logarithm of the Beltrami coefficient is a harmonic function, then the target surface can be equipped with a metric, conformal to the original one, and the solution of the Beltrami equation is a harmonic map. Therefore, solving a certain Beltrami equation is equivalent to solving the harmonic map problem. 
Harmonic maps to a constant curvature surface are therefore classified by the classification of the solutions of the elliptic sinh-Gordon equation. 

The general problem of solving the sinh-Gordon equation, and then the corresponding Beltrami equation, is still open. Different well known harmonic maps to the hyperbolic plane are proved to be related to the one-soliton solutions of the elliptic sinh-Gordon equation. Moreover, an example is proposed which does not belong to the one-soliton solution of the elliptic sinh-Gordon equation.  Solutions are calculated for the constant curvature case in a unified way, for positive, negative and zero curvature of the target surface.
\end{abstract}

\maketitle

\section{Introduction and Statement of the Results}\label{Introduction}

The aim of this article is to develop a method to construct harmonic diffeomorphisms, with nonvanishing Hopf differential, between Riemann surfaces $M$ and $N$. The case when $N$ is of constant curvature is studied in more detail. The method to obtain a harmonic map is summarized as follows: 

a) Find a solution of the elliptic sinh-Gordon equation.

b) Solve the Beltrami equation.  

c) Describe the metric on $N$ of constant curvature. 

There are  several examples of harmonic diffeomorphisms with nonvanishing Hopf differential, see for example \cite{C-T, Li-Ta2,S-T-W, Wolf1, Wolf2, X-A}. Recall that a minimal surface in $\mathbb{H}^2\times \mathbb{R}$ projects to a harmonic map to $\mathbb{H}^2$.  There are many examples of such harmonic maps obtained by methods similar to the study of the sinh-Gordon equation. Using the proposed method and the elliptic functions, we can find a family of harmonic maps to constant curvature spaces. This family includes some of the above examples and generalizes them. 

The study of harmonic diffeomorphisms between two Riemann surfaces is central to the theory of harmonic maps (see for example \cite{S-Y, Hitchin} and note that a twisted harmonic map, is a harmonic map in a local sense and that each harmonic map induces a solution of Hitchin self-duality equations). The case that has been studied the most is when the surfaces are of constant curvature (see for example \cite{A-R,Kal,Li-Ta2,M2,Minsky, P-H,S-T-W,  Wolf1, Wolf2, Wolf3} and the references therein). The geometric behavior of harmonic maps between hyperbolic surfaces has been studied in \cite{Han, Minsky, Wolf1, Wolf2, Wolf3}. The preparation of this article was motivated by the work in \cite{M}, that proves a conjecture of R. Schoen on quasiconformal harmonic diffeomorphisms between hyperbolic spaces.

In this paper, it is shown that the study of harmonic maps reduces to the study of the Beltrami equation of a certain type: the imaginary part of the logarithm of the Beltrami coefficient, coincides with the imaginary part of the logarithm of the Hopf differential, it is therefore a harmonic function. In addition, the real part of the logarithm of the Beltrami coefficient satisfies an elliptic nonlinear differential equation, which in the case of constant curvature is an elliptic sinh-Gordon equation. This last result reduces the harmonic map equations to a special case of the Beltrami equation that has been already extensively studied, see for instance \cite{G-R-S-Y}. 

The aforementioned harmonic map equations also feature for example in \cite{Han, H-T-T-W, Minsky, S-Y, Wolf1, Wolf2, Wolf3}. Note finally that there is a similar analysis of the Wang equation, which is applied in the study of affine spheres, see for example \cite{Loftin, DW} .

Furthermore, one could apply any property or calculation on a harmonic map to the associated Beltrami equation and vice versa. This result connects two previously extensively studied aspects of mathematics, i.e. harmonic maps and the Beltrami equation.

Harmonic maps to surfaces of constant curvature are closely related to the elliptic sinh-Gordon equation. The sinh-Gordon  and sine-Gordon equations have many applications and they have both been the subject of extensive study (see for example \cite{F-P,H}). Note that the sinh-Gordon equation was also crucial to the breakthrough work \cite{Wente} on the Wente torus. A close relation to the theory of constant mean curvature surfaces has been already known, \cite{Joa,K}. 

There are several model solutions of the sinh-Gordon equation, such as the soliton solutions, solutions by using separation of variables etc. All these models in their turn, imply models of solutions of the associated Beltrami equation and finally of the corresponding harmonic map problem.  Using one-soliton solutions of the elliptic sinh-Gordon equation one can construct examples of harmonic maps to constant curvature spaces. Elliptic functions arise in course of this construction. With this approach, one can study positive, negative and zero constant curvature surfaces in a unified way, and find concise formulas that provide new examples of harmonic maps. As an application, one can recover some of the examples of harmonic maps presented in the articles \cite{Li-Ta2, S-T-W,Wolf1, Wolf2}, by proving that they correspond to the one-soliton solution of the sinh-Gordon equation.

Finally, a B{\"a}cklund transform arises, which provides a connection between the solutions of a certain elliptic sinh-Gordon and an elliptic sine-Gordon equation. In fact, this result provides solutions of the first equation by the known solutions of the second and vice versa. As an application, a new harmonic map is constructed, that does not correspond to one-soliton solutions of the sinh-Gordon equation.

Throughout this article, we assume that the Hopf differntial of the harmonic map does not vanish. The main results in this article could be summarized in the following Theorems.

  \begin{theorem}\label{Prop:Beltrami_to_harmonic}
If the diffeomorphism $u: M \to N$ satisfies the Beltrami equation \[\frac{u_{\bar{z}}}{u_{{z}}}= e^{ -2 \omega(z,\bar{z})+i\phi(z,\bar{z})},\]
  where $\phi_{{z}\bar{z}}=0,$ then $N$ can be equipped with a conformal metric such that $u$ is a harmonic map and the curvature of $N$ is
   \[
  \tilde{K}_{N}=-\frac{2\omega_{z\bar{z}}}
  {\sinh{2\omega}}e^{\psi},
  \]
   where $\psi$ is the conjugate harmonic function to $\phi$.
 \end{theorem}

\begin{theorem}
	\label{theorem_main_equivalence}
 A necessary and sufficient condition for a harmonic diffeomorphism $u: M \to N$ to exist, is that the imaginary part of the Beltrami coefficient is a harmonic function.   
\end{theorem}

The next corollary follows from the proof of Theorem 2. 

\begin{corollary}\label{main_corollary}     
  Let $\omega$ be a solution of the elliptic sinh-Gordon equation 
 \begin{equation*}
       \omega_{\zeta\bar{\zeta}}= - \frac{K_N}{2}  \sinh 2 \omega,
    \end{equation*}
    and let $u$ be a solution of the Beltrami equation 
    \begin{equation}\label{Beltrami eqn}
\frac{u_{\bar{\zeta}}}{u_{{\zeta}}}= e^{ -2 \omega(\zeta,\bar{\zeta})}.
\end{equation}
Then, a harmonic map $U=U(\zeta,\bar{\zeta})$ that satisfies (\ref{Beltrami eqn})
 can be written as
   \begin{equation*}
 U=f\left(  u(\zeta,\bar{\zeta})\right), 
\end{equation*}
where $f(z)$  is holomorphic and the metric on $N$ is of constant curvature $K_{N}$. Theorem \ref{theorem_main_equivalence} implies that there is a classification of  harmonic diffeomorphisms via the classification of the  solutions of the sinh-Gordon equation.             
\end{corollary}       

Let us now present an outline of the article. In Section \ref{sec:Background}, the necessary formulas are introduced. Next, Section \ref{sec:Proofs} contains the proof of Theorem \ref{theorem_main_equivalence}. In Section \ref{Beltrami}, the solution of the Beltrami equation is discussed and in Section \ref{sec:Constant}, the constant curvature case is studied by using one-soliton solutions of the sinh-Gordon equation. In Section \ref{sec:S-T-W}, it is shown that the solutions given in \cite{Li-Ta2,S-T-W,Wolf1,Wolf2}, correspond to the one-soliton solutions of the sinh-Gordon equation and explicit formulas are given. The results of this paper can be extended to the positive curvature case by using the explicit formulas of Section \ref{sec:Constant}, since those are given in a unified way for positive, negative and zero curvature. The B\"acklund transform is discussed in Section \ref{sec:Backlund}. Next, Section \ref{sec:prespectives} contains some perspectives for future research. Finally, there is an Appendix about elliptic functions.

\section{Preliminaries}\label{sec:Background}

\subsection{Isothermal Coordinates}

Let $u: M \rightarrow N$ be a map between Riemann surfaces $(M,g), (N,h).$ The map $u$ is locally represented by $u=u(z)=R+iS.$ The standard notation is  
\[
\partial_{z}=\frac{1}{2}(\partial_{x}-i\partial_{y}), \quad \partial_{\bar{z}}=\frac{1}{2}(\partial_{x}+i\partial_{y}), \quad z=x+iy.
\]

It is a known fact the existence of isothermal coordinates on an arbitrary surface with a real analytic metric (see \cite[Section 8, p. 396]{J}).  Consider an isothermal coordinate system $(x,y)$ on $M$ such that
 \begin{equation*}
g=e^{f(x,y)}(dx^2 + dy^2)=e^{f(z,\bar{z})} dz d{\bar{z}}=e^{f(z,\bar{z})}  |dz|^2,
\end{equation*} 
where $z=x+iy.$
Consider an isothermal coordinate system $(R,S)$ on $N$ such that
\begin{equation*}
h=e^{F(R,S)}(dR^2 + dS^2)=e^{F(u,\bar{u})}  du d{\bar{u}}=e^{F(u,\bar{u})}|du|^2,
\end{equation*}
where $u=R+iS.$

Note that the Gauss curvature on the target is given by 
\begin{equation}\label{curvature}
K_N=K_N(u,\bar{u})
=-\dfrac{1}{2} \mathrm{e}^{-F} \Delta F=-2F_{u\bar{u}} e^{-F}. 
\end{equation} 
 
 \subsection{Harmonic Maps and the Beltrami Equation}

In the case of isothermal coodinates (see \cite[Section 8, p. 397]{J}), a map $u$ is harmonic if it satisfies 
\begin{equation}\label{eq:harmonic_map}
u_{z\bar{z}} + F_{u}(u,\bar{u}) u_{z} u_{\bar{z}}=0.
\end{equation}
Note that this equation only depends on the complex structure of $N$ and not on the metric $g$ of $M.$

Denote by
 \begin{equation*}
\| u_z\|^2=e^{{F(u,\bar{u})}}e^{-{f(z,\bar{z})}} | u_z|^2 
 ,\quad
   \| u_{\bar{z}} \|^2= e^{{F(u,\bar{u})}}e^{-{f(z,\bar{z})}} | u_{\bar{z}} |^2 
   \end{equation*}
the norms of the (1,0)-part and (0,1)-part of $du,$ respectively.

The Jacobian of $u$ is defined by 
\begin{equation}\label{eq:notation_Jacobian}
\begin{split}
J(u) =& \| u_z\|^2 -\| u_{\bar{z}} \|^2 \\
= &e^{{F(u,\bar{u})}}e^{-{f(z,\bar{z})}}\left( |u_z |^2 -| u_{\bar{z}} |^2\right).
\end{split}
\end{equation}
Note that if $u\colon M \rightarrow N$ is a diffeomorphism, then the Jacobian 
 is nowhere vanishing. 

The Hopf differential of $u$ is defined by
\begin{equation}\label{eq:Hopf_theorem}
\Lambda(z) dz^2=\left( e^{F(u,\bar{u})}u_{z} \bar{u}_{z}\right)dz^2.
\end{equation}
It is a well known result  \cite[Section 8, p. 399]{J} the following proposition.
\begin{proposition}\label{prop:Hopf}
A necessary and sufficient condition for  $u$ to be a harmonic map, is that the Hopf differential of $u$ is holomorphic, i.e.
\[
e^{F(u,\bar{u})} u_z    \bar{u}_z= e^{-\lambda(z)},
\]
where $\lambda(z)$ is a holomorphic function.\end{proposition}

Consider the Beltrami coefficient
\begin{equation*}
\mu(z,\bar{z})=\frac{u_{\bar{z}}}{u_{{z}}}.
\end{equation*}
Then, the following relations are valid:
\begin{equation*}
du= u_{{z}} dz+ u_{\bar{z}} d \bar{z} = u_{{z}} \left(  dz + \mu(z,\bar{z}) d\bar{z} \right),
\end{equation*}
where 
\begin{equation*}
\frac{u_{\bar{z}}}{u_{{z}}}=\mu(z,\bar{z})= e^{ -2 \omega(z,\bar{z})+ i \phi(z) }.
\end{equation*}
This is the well known Beltrami equation. Note that, in general, the Beltrami coefficient $\mu(z,\bar{z})$ is a complex function.

\section{Proof of the Theorems}\label{sec:Proofs}

 The following formulation will be used in the text:
if  $u$ is a harmonic map, then 
\begin{equation}\label{eq:Hopf}
e^{F(u,\bar{u})}  u_z   \bar{u}_z= e^{-\lambda(z)}  
\mbox{ and }  
 e^{F(u,\bar{u})} u_{\bar{z}} \bar{u}_{\bar{z}} =e^{-\overline{\lambda(z)}},
\end{equation}
where $\lambda(z)$ is a holomorphic function.



  Motivated by (\ref{eq:Hopf}), we set
   \begin{equation}\label{eq:def_w}
   e^v = e^{ \frac{F+\lambda}{2}} u_z \mbox{ and  } v=\omega+i \theta .
   \end{equation}
     Then 
  \begin{align}\label{eq:exp_v}
  e^{v} = e^{ \frac{F}{2}}\, e^{\frac{\lambda}{2}} {u}_z= 
 \,\dfrac{1}{ e^{ \frac{F}{2}} e^{\frac{\lambda}{2}} {\bar{u}}_z}
 & \mbox{ , } & 
   e^{-v }= e^{ \frac{F}{2}}\, e^{\frac{\lambda}{2}} \bar{u}_z= 
 \,\dfrac{1}{ e^{ \frac{F}{2}} e^{\frac{\lambda}{2}} {u}_z}
  \\
  \label{eq:exp_vbar}
  e^{\bar{v}} = e^{ \frac{F}{2}}\, e^{\frac{\bar{\lambda}}{2}} {\bar{u}}_{\bar{z}}= 
 \,\dfrac{1}{ e^{ \frac{F}{2}} e^{\frac{\bar{\lambda}}{2}} {u}_{\bar{z}}}
 & \mbox{ , } & 
   e^{-\bar{v} }= e^{ \frac{F}{2}}\, e^{\frac{\bar{\lambda}}{2}} {u}_{\bar{z}}= 
 \,\dfrac{1}{e^{ \frac{F}{2}} e^{\frac{\bar{\lambda}}{2}} {\bar{u}}_{\bar{z}}}
 \end{align}
 and
 \begin{equation*}
 e^{F+\lambda} u_z  \bar{u}_z=1,   \quad
   e^{F+\bar{\lambda}} \bar{u}_{\bar{z}}  \bar{u}_{\bar{z}}= 1 .
\end{equation*}
The above equations imply
\begin{equation}\label{eq:relations_w}
e^{-(v+\bar{v})}=e^{-2 \omega}= 
\dfrac{ e^{\frac{\bar{\lambda}}{2}} {u}_{\bar{z}} } {e^{\frac{{\lambda}}{2}} {u}_{{z}} }   
=
\dfrac{ e^{\frac{{\lambda}}{2}} {\bar{u}}_{{z}} } {e^{\frac{{\bar{\lambda}}}{2}} {\bar{u}}_{\bar{z}} }   ,
\end{equation}
therefore the harmonic map satisfies the Beltrami equation:
 \begin{equation}\label{eq:Beltrami_w}
\dfrac{   {u}_{\bar{z}} }{   {u}_{z} }= e^{ - 2 \omega+  i \, \mathrm{ Im}\, \lambda(z) }= \mu(z,\bar{z}).
 \end{equation}
 Note that the imaginary part of the logarithm of the Beltrami coefficient  $\mu$ is a harmonic function, since it is the imaginary part of a conformal function. 
 
 Equations (\ref{eq:def_w}),   (\ref{eq:exp_v})  and (\ref{eq:exp_vbar})  imply that
\begin{equation*}
e^{2 \omega} = e^{v+\bar{v}}= e^{F+\mathrm{Re}\, \lambda} | u_z|^2,  \quad 
e^{-2 \omega} = e^{-(v+\bar{v})}= e^{F+\mathrm{Re}\, \lambda} | u_{\bar{z}}|^2
\end{equation*}   
and
\begin{equation}\label{eq:sinhw1}
\sinh(2 \omega)=  \dfrac{ e^{F+\mathrm{Re}\, \lambda}}{2}  \left(
|u_z|^2 -|u_{\bar{z}}|^2\right),   \;
 \cosh(2 \omega)=  \dfrac{ e^{F+\mathrm{Re}\, \lambda}}{2}  \left(
|u_z|^2+|u_{\bar{z}}|^2\right).
\end{equation}       
 As a consequence of (\ref{eq:exp_v}), it follows that
 \begin{equation*}
 e^{2 v}= \frac{u_z}{{\bar{u}}_z}, \text{ thus } 2v = 2 \omega +2 i \theta =\log u_z -  \log{\bar{u}_{{z}}}.
 \end{equation*}      
 After some elementary  calculations (see Section \ref{sec:Backlund} for more details) and taking into consideration (\ref{eq:harmonic_map}), it follows that
 \begin{equation}\label{eq:v_zbarz}
    \begin{split}
        2 v_{ z \bar{z}} =& 2 \omega_  { z \bar{z}} +2 i \theta_{ z \bar{z}}  =\\
        =&F_{ u \bar{u}} \left( |u_z|^2 - | u_{\bar{z}}|^2\right) +  \left( F_u^2-F_{uu}\right) u_z u_{\bar{z}} -  \left( F_{\bar{u}}^2-F_{\bar{u} \bar{u}}\right) {\bar{u}}_z {\bar{u}}_{\bar{z}} .
        \end{split}
\end{equation}
Therefore, from (\ref{curvature}), (\ref{eq:sinhw1}) and (\ref{eq:v_zbarz}), we have
\[
2\omega_{z\bar{z}} = F_{u\bar{u}} \left(|u_z|^2- |u_{\bar z}|^2 \right) 
=2 F_{u \bar{u}} e^{-F-\mathrm{Re\lambda} }\sinh 2 \omega= - K_N e^{-\mathrm{Re} \lambda} \sinh 2 \omega,
\]
    where $K_N$ the curvature of the surface $N$.

The following is a fundamental result in the theory of harmonic maps that can be found in \cite{Minsky, Wolf1, Wolf3}. The proof has been included above for clarity reasons. 
 \begin{proposition}[Minsky\cite{Minsky}, Wolf\cite{Wolf1}]
 	\label{Prop:harmonic_to_Beltrami}
 If $u$ is a harmonic map then it satisfies the Beltrami equation
  \begin{equation}\label{Beltrami2}
   \frac{u_{\bar{z}}}{u_{{z}}}= e^{ -2 \omega(z,\bar{z})+i\phi},
  \end{equation}
  and $\phi$ is a harmonic function i.e. $\phi_{{z}\bar{z}}=0$.
   
Furthermore, if $\psi$ is the conjugate harmonic function to $\phi$, then 
   \[K_{N}=-\frac{2\omega_{z\bar{z}}}{\sinh{2\omega}}e^{\psi},
  \]
  where $K_N$ is the curvature of the target manifold N.
 \end{proposition}


 Suppose now that the diffeomorphism $u$ is a solution of this Beltrami equation (\ref{Beltrami2}), such that
 \begin{equation*}
 \mathrm{Re}\,\log(\mu(z,\bar{z}))= -2 \omega(z,\bar{z}),\;
 \mathrm{Im}\,\log(\mu(z,\bar{z}))=  \phi(z,\bar{z})\; \mbox{ and }  \phi_{z\bar{z}  }=0.
 \end{equation*}
 Then,  there exists a holomorphic function $\lambda(z)=\psi(z)+ i\phi(z),$ where $\psi(z)$ is the harmonic conjugate of $\phi(z)$.
 
 We define the function
  \begin{equation*}
 	\tilde{F}(u,\bar{u})=\log\frac{e^{-\lambda(z)}}{u_z  \bar{u}_z}.
\end{equation*}	
From the Beltrami equation, it follows that
 \begin{equation}
 \label{eq:tildeF}
 e^{\tilde{F}(u,\bar{u})}= \dfrac{e^{-\psi}}{|u_{z}|^2}> 0.
 \end{equation}
Given the imaginary part of the  logarithm of the Beltrami equation, which is a harmonic function, and the fact that the map $u$ is a diffeomorphism, we can calculate $\tilde{F}$ as a function of $u$ and $\bar{u}$.
 
  Let us equip $N$ with the conformal metric 
  \begin{equation}\label{eq:metric}
  \tilde{h}=e^{\tilde{F}}ds^2.
  \end{equation} 
Then, 
\[
 e^{\tilde{F}} u_z  \bar{u}_z=e^{-\lambda(z)},
\] 
 thus
 \begin{equation*}
 u_{z\bar{z}} + \tilde{F}_u (u,\bar{u})  u_{z}   u _{\bar{z}}=0 ,
 \end{equation*}
 i.e. $u$ is a harmonic map. Then, from Proposition \ref{Prop:harmonic_to_Beltrami} follows that the corresponding curvature is 
 \begin{equation}\label{star}
 \tilde{K}_{N}=
 -\frac{2\omega_{z\bar{z}}}{\sinh{2\omega}}e^{\psi}.
\end{equation} 
Therefore Theorem \ref{Prop:Beltrami_to_harmonic} is valid. Note that Theorem \ref{Prop:Beltrami_to_harmonic} is the converse of Proposition \ref{Prop:harmonic_to_Beltrami}, when the map is a diffeomorphism. These two results could be summarized to Theorem \ref{theorem_main_equivalence}.

Note that when
    $\lambda(z)=\psi(z)+ i\phi(z)$ is zero and $N$ is a surface with constant curvature, then the equation (\ref{star}) reduces to
    \begin{equation}
        \omega_{z\bar{z}}= - \frac{K_N}{2} \sinh 2 \omega,
    \end{equation}       
which  is the well known elliptic sinh-Gordon differential equation, with many applications in physics. There is an extensive bibliography on the solutions of this equation, see for example \cite{H} and the references therein.
    In the next sections, we shall discuss the construction of harmonic maps between surfaces by using the solution of the elliptic sinh-Gordon equation.

     In the domain surface $M$ of a harmonic map, one can choose a specific coordinate system in order to considerably facilitate the calculations. This \textit{specific}  system is defined by the conformal transformation
\begin{equation}\label{eq:specific}
\zeta=\xi+ i \eta= \int e^{-\lambda(z)/2} \, dz,
     \end{equation}     
 where     $\lambda(z)$ is the holomorphic function given by the Hopf differential (\ref{eq:Hopf_theorem}) and it is related to the imaginary part of the logarithm of  the Beltrami coefficient $\mu,$ see equation (\ref{eq:Beltrami_w}). In this specific system the equations in Section \ref{sec:Proofs} could be simplified by substituting $\lambda(z)=0$.

  Then,
  \begin{equation}\label{eq:holomorphic_map_specific}
u_{\zeta\bar{\zeta}} + F(u,\bar{u})_{u} u_{\zeta} u_{\bar{\zeta}}=0 , \quad 
e^{F(u,\bar{u})}  u_\zeta   \bar{u}_\zeta= 1.
  \end{equation}
 The corresponding Beltrami equation is given by 
 \begin{equation}\label{eq:Beltrami_Equation_specific1}
 \frac{u_{\bar{\zeta}}}{u_{{{\zeta}}}}= e^{ -2 \omega(\zeta,\bar\zeta)} ,
 \end{equation}
 where the function $\omega(\zeta, \bar{\zeta})$ satifies the elliptic sinh-Gordon equation
 \begin{equation}\label{eq:sinh-Gordon_specific1}
    \omega_{\zeta\bar{\zeta}}= - \frac{K_N}{2}    \sinh 2 \omega.
 \end{equation}    
 The above equations are of special interest in the case of constant sectional curvature of the target manifold $N$. In this case, given a function that satisfies the elliptic sinh-Gordon equation (\ref{eq:sinh-Gordon_specific1}), one has to calculate a solution of the Beltrami equation (\ref{eq:Beltrami_Equation_specific1}). By the conformal change of coordinates  (\ref{eq:specific}) of this solution, one can calculate the solution in the original coordinates $z=x+i y$ and the general solution of the problem is a holomorphic function of this solution.  The above is a strategy to solve the harmonic map problem from a domain surface $M$ to a target surface $N$ of constant curvature. Thus Corollary \ref{main_corollary} is valid and it implies that there is a classification of  harmonic diffeomorphisms via the classification of the  solutions of the sinh-Gordon equation.

\section{Solution of the Beltrami Equation}\label{Beltrami}

Let us consider a solution $\omega$ of the elliptic sinh-Gordon equation
equation (\ref{eq:sinh-Gordon_specific1}), then the Beltrami equation (\ref{eq:Beltrami_Equation_specific1}) can be written as a first order P.D.E.
\begin{equation}\label{eq:Compatibility}
 e^{  \omega(\zeta,\bar{\zeta})} u_{\bar{\zeta}} - e^{ - \omega(\zeta,\bar{\zeta})}  u_{\zeta}=0.
\end{equation}
The relations $\zeta=\xi+i\eta$  and $u=R(\xi,\eta)+ i S(\xi,\eta)$ 
and the fact that $\omega(\zeta, {\bar{\zeta}}  )$  is a real function, imply that the following system  of first order  PDE's holds true:
\begin{align}
\sinh\, \omega\, R_\xi - \cosh\, \omega\, S_\eta =0 \label{eq:system_1}
\\
\cosh\, \omega\, R_\eta +\sinh\, \omega\, S_\xi=0  . \label{eq:system_2}
\end{align}
By multiplying the equation (\ref{eq:system_1}) by $R_\eta$ and  equation (\ref{eq:system_2}) by $S_\eta$ the following Proposition holds true.

\begin{proposition}\label{prop:orthogonality_Beltrami_solutions}
The solution of the Beltrami equation with a real Beltrami coefficient $\mu(\zeta, \bar{\zeta})$  corresponds to a harmonic map which preserves the orthogonality of the local coordinate systems, i.e.
\begin{equation*}
R_\xi R_\eta \, + \, S_\xi S_\eta=0.
\end{equation*}
\end{proposition}
The system (\ref{eq:system_1}, \ref{eq:system_2})  can be separated as two second order O.D.E.'s
\begin{align}
\left(  \tanh\,\omega\, R_\xi \right)_\xi +\left(  \coth\,\omega\, R_\eta \right)_\eta =0
\label{eq:2_system_1}
\\
\left(  \tanh\,\omega\, S_\xi \right)_\xi +\left(  \coth\,\omega\, S_\eta \right)_\eta =0.
\label{eq:2_system_2}
\end{align}
Indeed, from (\ref{eq:system_1}) we get $ S_\eta= \tanh\, \omega\, R_\xi $, from (\ref{eq:system_2}) we get $ S_\xi=- \coth\, \omega\, R_\eta,$ and from equality $\left(S_{\xi}\right)_{\eta}= \left(S_{\eta}\right)_{\xi},$ the equation (\ref{eq:2_system_1}) follows. The proof of equation (\ref{eq:2_system_2}) is similar, thus omitted. 
Then, equation (\ref{eq:2_system_1}) can be written as
\begin{equation}\label{eq:2nd_order}
\tanh\, \omega\,   R_{\xi\xi}+  \coth\, \omega\,   R_{\eta\eta}
+ \dfrac{ \omega_\xi\,   R_\xi}{ \cosh^2\,\omega}- \dfrac{ \omega_\eta\,   R_\eta}{ \sinh^2\,\omega}=0.
\end{equation}
This is an elliptic  second order P.D.E. Given a solution $\omega(\xi,\eta)$ of the elliptic sinh-Gordon equation, that is 
\begin{equation}
    \omega_{\xi\xi}+  \omega_{\eta\eta}= - 2K_N (\xi,\eta  )  \sinh 2 \omega,
 \end{equation}    
then one can solve equation (\ref{eq:2nd_order}), by using standard methods, see for example \cite[p. 72]{P-R}. In other words, one has to solve the equation of characteristics
\begin{equation}\label{eq:character}
\dfrac{d \eta}{d\xi}= i \epsilon \coth \omega(\xi,\eta), \; \epsilon=\pm 1.
\end{equation}
The solution of the above equation is of the form
\[
\Phi(\xi,\eta)= \mathrm{constant}.
\]
Then, the solutions of (\ref{eq:Compatibility}) are given by
 \[
 u= \mathrm{Re} \Phi+ i  \mathrm{Im }\Phi=R(\xi,\eta)+ i S(\xi,\eta) .
 \]

Thus, we once again emphasize the strategy to find solutions of the harmonic map problem: it is enough to take a solution of the elliptic sinh-Gordon equation and find a solution of the Beltrami equation (\ref{eq:Compatibility}). 



\section{Constant curvature spaces, one-soliton solution}\label{sec:Constant}

In this section we consider the case when $N$ is of constant curvature $K_{N}$. In the specific coordinates (\ref{eq:specific}), the elliptic sinh-Gordon equation is 
\[
\omega_{\zeta \bar{\zeta}}=-\frac{K_N}{2} \sinh 2 \omega, \mbox{ where } K_N=\pm 1, 0.
\]
A general solution for this equation is not known. There are only partial solutions, and of particular interest is the so called one-soliton solution, defined by:
\[
\omega=\omega( \gamma\eta-\delta\xi  ),   \;  \gamma=\rho \cos \tau  \mbox{  and  }  \delta=\rho \sin 
\tau .
\]
In order to simplify the calculations, a new systen of coordinates $Z=X+i Y$  is introduced:
\begin{equation}\label{defn rho and tau}
\begin{array}{c}
X= \xi \rho \cos\tau+ \eta \rho \sin \tau
\\
Y= -\xi \rho \sin\tau+ \eta \rho \cos \tau
\end{array}
\mbox{ or }
\begin{array}{c}
Z= \rho \mathrm{ e }^{- i\,\tau} \zeta
,\;
\bar{Z}= \rho \mathrm{ e }^{i\, \tau} \bar{\zeta}.
\end{array}
\end{equation} 
In these coordinates, $\omega=\omega(Y)$ and the elliptic sinh-Gordon equation is written
\begin{equation}\label{eq:mod_Sinh_G}
\dfrac{  d^2 \omega}{dY^2}=-\dfrac{ 2K_N}{\rho^2} \sinh 2\omega  \; \mbox{  or  }  \;
\left(\dfrac{  d \omega}{dY}\right)^2=
C\left(
1 
 - \dfrac{ 4 K_N} {C\rho^2}
	 \sinh^2 \,\omega
\right).
\end{equation}
Equivalently
\begin{equation}\label{eq:der_omega_equation}
\left(\dfrac{\omega '(Y)}{\sqrt{C}}                                                   
\right)^2+ (m-1)  
\sinh^2 \,\omega(Y)=1,
\end{equation}
where 
\begin{equation}\label{constants}
\omega_0=\omega(Y_0), \; \omega '_0=\omega '(Y_0), 
\end{equation}
and
\[
C= \left(\omega '_0                                                   
\right)^2+ \dfrac{4K_N} {\rho^2}
\sinh^2 \,\omega_0, 
\; m= 1+  \dfrac{4K_N} {C \rho^2}.
\]
The parameter $Y_0$ corresponds to the choice of the initial conditions. The solution of the  equation  (\ref{eq:der_omega_equation}) can be calculated by using the Jacobi elliptic functions, discussed in the Appendix. More precisely, it follows that
\begin{equation}\label{omega}
\dfrac{\omega '(Y)}{\sqrt{C}}                                                   
= cd(\sqrt{Cm}(Y-Y_0)+v_0 \vert \frac{1}{m}),
\end{equation}
where 
\begin{equation*}
v_0 = sd^{-1}(\sqrt{m}\sinh{\omega_0}\vert \frac{1}{m}).
\end{equation*}
Moreover,
\begin{equation}\label{tanhomega}
\tanh{\omega}=\frac{1 }{ \sqrt{m}} sn(    \sqrt{Cm}(Y-Y_0)+v_0 \vert \frac{1}{m}).
\end{equation}
In the coordinates $Z=X+i Y$, the equation (\ref{eq:Compatibility}) is written as
\begin{equation}\label{eq:Compatibility_Z}
e^{  \omega(Y)+ i \tau} u_{\bar{Z}}-e^{ - \omega(Y)- i \tau}  u_Z=0.
\end{equation}
This equation is equivalent to the equation
\begin{equation}\label{eq:2Compatibility_Z}
i u_{Y} \coth{(\omega (Y)+ i \tau)} + u_{X}=0\end{equation}
One obvious solution of the above equation is given by
\begin{equation*}\label{eq:solutionY}
u=  R+ i S=R_0+ i S_0+
\alpha\left({
X-X_0 + i   \int_{Y0}^{Y} \tanh( \omega(t) + i  \tau  ) dt
}\right),
\end{equation*}
since 
\[
u_{X}=\alpha \text{ and } u_{Y}=i\alpha \tanh{(\omega (Y)+ i \tau)}.
\]

Note that
\begin{equation}
\label{eq:solution_semi}
R-R_0=\alpha\left({X-X_0 +  \mathrm{Re} \left( i   \int_{Y0}^{Y} \tanh( \omega(t) + i  \tau  ) dt\right)}\right),
\end{equation}
\begin{equation}
\label{eq:solution_semi_1}
S-S_0= \alpha\left({\mathrm{Im} \left( i   \int_{Y0}^{Y} \tanh( \omega(t) + i  \tau  ) dt\right)}\right).
\end{equation}
Elementary calculations 
give that
\begin{equation}\label{tanh}
 i   \tanh( \omega(t) + i  \tau  ) =
 \dfrac{-\tan (\tau )+i   \tanh\omega}{1+i   \tan (\tau ) \tanh
  \omega}.
\end{equation}
Differentiating (\ref{eq:solution_semi}) and applying (\ref{tanh}), we obtain
\begin{equation}\label{eq:RY}
\begin{array}{rl}
\dfrac{\partial R}{\partial Y}&=\alpha \, \mathrm{Re} \,\left[ i  \tanh( \omega + i     \tau  )\right] = - \dfrac{ \alpha (1-\tanh^2{\omega})  \tan \tau } {  1+ \tan^2 \tau \tanh^2\omega }=\\
&=-	\dfrac{(m-1)  \sin\, \tau\,  \cos\,\tau  }{M}\dfrac{\alpha}{ 1-\left(  \dfrac{\omega'(Y)}{\sqrt{CM}}\right)^2  } ,
 \end{array}
\end{equation}
\\
\\
 { where }
 \[
 M= 1+\dfrac{4K_N} {C\rho^2}\cos^2 \tau= \dfrac{m+\tan^2 \tau}{  1+\tan^2\tau   }.
\]
Also, differentiating (\ref{eq:solution_semi_1}), we obtain
\begin{equation}\label{eq:SY}
\begin{array}{rl}
\dfrac{\partial S}{\partial Y}&=\alpha \mathrm{Im}\left[ i   \tanh( \omega  + i  \tau  )\right]=    \dfrac{ \alpha    (1+\tan^2{\tau}) \tanh\omega } {  1+ \tan^2 \tau \tanh^2\omega }=\\
&=
 \alpha  \dfrac{C }{M   }
 \, \dfrac{   \omega''(Y)  }{ 1-\left(  \dfrac{\omega'(Y)}{\sqrt{CM}}\right)^2  } .
\end{array}
\end{equation}
We find
\begin{equation*}\label{eq:eF}
\begin{array}{rl}
\Phi=e^{F}=& \dfrac{1}{u_{\zeta} \bar{u}_{{\zeta}}}=
\dfrac{ 4}{ \alpha^2  \rho^2 }\left({  \cos^2 \tau \cosh^2 \omega+ \sin^2 \tau \sinh^2 \omega }\right).
\end{array}
\end{equation*}
One should note an interesting relation:
\begin{equation}\label{eq:relation_eFRY}
 \dfrac{\partial R}{\partial Y}=-\dfrac{4  \sin\, \tau\, \cos\, \tau}{\alpha \rho^2}  \dfrac{1}{\Phi}.
\end{equation}
Another interesting relation, which shall be applied in Section \ref{sec:S-T-W}, is the following one:
\begin{equation}\label{eq:rel}
\left(\dfrac{\partial S}{\partial Y}\right)^2=
-\dfrac{16 {\tan^2 \,\tau }}{\left(\tan^2 \, \tau +1\right)^2
   \alpha ^2 \rho ^4}\dfrac{1}{\Phi^2}
   +\dfrac{4 \left({\tan^2\,\tau
   }-1\right)}{\left({\tan\,\tau  }^2+1\right) \rho ^2 \Phi
   }+\alpha ^2.
\end{equation}
Therefore (see the  detailed explanations and  calculations in the Appendix),
\begin{equation}\label{eq:R_Pi}
\begin{array}{l}
R-R_0=
\alpha (X-X_0) -  \\
- \epsilon \sqrt{ C m} 	\dfrac{(m-1) \,\alpha \, \sin\, \tau\,  \cos\,\tau  }{M}
\left(
\Pi\left( \dfrac{1}{M}; \,v +K \,\vert\,\dfrac{1}{m}  \right)-
\Pi\left(
 \dfrac{1}{M}; \,v_0 +K \,\vert\,\dfrac{1}{m} 
\right)
\right)
\end{array}
\end{equation}
where
\[
v=\epsilon \sqrt{Cm}(Y-Y_0)+v_0 \;  \mbox{and}   \;  K=K(\dfrac{1}{m})=
\int_{0}^{\frac{\pi}{2}} 
\dfrac{d\theta}{\sqrt{ 1 -\dfrac{1}{m} \sin^2
 \theta}}
\]
where $\Pi(n,x\vert m)$ is the elliptic integral of the third kind. Also,
\begin{equation}\label{eq:S_fin}
S-S_0=
\dfrac{\alpha }{\sqrt{C M}}\left(
{
\mathrm{arctanh}\dfrac{\omega'(Y)}{\sqrt{C M} }-\mathrm{arctanh}\dfrac{\omega'_0}{\sqrt{C M}} 
}
\right)
\end{equation}
and
\begin{equation*}
e^F= \dfrac{ 4   M} {(m-1) \alpha^2  \rho^2 \cosh^2{\Sigma}} , 
\end{equation*}
where
\begin{equation}
\Sigma= \frac{\sqrt{C M}}{\alpha} (S-S_0) +  \mathrm{arctanh}\dfrac{\omega'_0}{\sqrt{C M}}.
\end{equation}
Note that the metric on $N$ is of constant curvature and that the results in Section \ref{sec:Constant} cover all cases of positive, negative and zero constant curvature in a unified formulation.


\section{Explicit solutions}\label{sec:S-T-W}

\subsection{The hyperbolic cylinder examples of Wolf \cite{Wolf1,Wolf2}}
This Section focuses on the explicit solutions of the harmonic map problem in the influential works \cite{Wolf1,Wolf2}. These results can be recovered by the results in Section \ref{sec:Constant}. 

Consider $M$ to be a hyperbolic cylinder with boundary that can be realized as the rectangle $[-1,1] \times [0,1]$ in the 2-plane with the Euclidean metric, identifying $[-1,1] \times \{0\}$ with $[-1,1] \times \{1\}$ to obtain a cylinder.  

Consider $N$ as the rectangle $$[-\cosh^{-1}t, \cosh^{-1}t] \times [0,1],\, t>0,$$
equipped with the metric $$ds^2 = du^2 + {\dfrac{\cosh^2 u}{t^2}} dv^2,$$ identifying $[-\arccos\mathrm{h} \,t, \arccos\mathrm{h} \,t] \times \{0\}$
  and $[-\arccos\mathrm{h} \,t, \arccos\mathrm{h} \,t] \times \{1\}$ to obtain a cylinder. 
  
  Then $N$ has constant curvature -1, and $M$ has the Euclidean conformal structure. 

In \cite{Wolf1}  a harmonic map is found, which takes the form $u(x,y) = u(x)+i y$ and $u$ satisfies the boundary value problem 
\begin{equation}
u'' (x) = \frac{1}{2t^2}\sinh{2u(x)}, \, u(0)=0, \, u(1)= \cosh^{-1}t.
\end{equation} 
It follows that 
\begin{equation}\label{wolf u}
\left( u' (x) t \right)^2 = \frac{1}{2}\cosh{2u(x)}+c_{0}=\sinh^{2}{u(x)}+\frac{1}{2}+c_{0},
\end{equation} 
for some constant $c_{0}$. This equation can be solved by using elliptic integrals, see the Appendix for the definition of the elliptic functions.

Under the change of variables $u=\sinh^{-1}\left(\tan{\frac{S}{t}}\right), \, v=R,$ the metric in the target becomes 
\[
h=e^{F}(dR^2+ dS^2)=\frac{1}{t^2\cos^2{\left(\frac{S}{t}\right)}}(dR^2+ dS^2).
\]
Then, the given harmonic map becomes
\begin{equation}
u(x,y)=y+ i t \arctan{\sinh{u(x)}},
\end{equation} 
where $u=u(x)$ is as above.  

Taking into account (\ref{wolf u}), it follows that the Hopf differential is given by
\begin{equation}
u_{z} \bar{u}_{z} e^{F}=\frac{1}{4}((u')^2 -\frac{\cosh^{2}u}{t^2})=\frac{1}{4t^2}(c_{0}  -1/2),
\end{equation}
see (\ref{wolf u}) above. The Beltrami coefficient is given by 
\begin{equation}
\frac{u_{\bar{z}}}{ u_{z}}=\frac{t u' + \cosh{u}}{t u' - \cosh{u}}=e^{-2 \omega},
\end{equation}
and $\omega$ satisfies the equation
\begin{equation}
\omega'' = \frac{4}{c_{0} t^2 -1/2}\sinh{2\omega}.
\end{equation}
Note that the choice of the initial conditions in (\ref{constants}) is as follows:
\begin{equation}
\omega_{0}=-\frac{1}{2}\log\left(   \frac{\sqrt{c_{0} t^2 +1/2}-1}{\sqrt{c_{0} t^2 +1/2}+1}\right),\, \omega'_{0}=0.
\end{equation}

Set
\begin{equation}
\zeta = \frac{2}{\sqrt{c_{0} t^2 -1/2}}z.
\end{equation}
Then, the solution becomes
\begin{equation}\label{Wolf1}
u(\xi,\eta)=\frac{2\eta}{\sqrt{c_{0} t^2 -1/2}}+ i t \arctan{\sinh{u(\frac{2\xi}{\sqrt{c_{0} t^2 -1/2}}}}).
\end{equation} 

The same harmonic map can be recovered by the method presented earlier. More precisely, consider \[\alpha=-\frac{2t}{\sqrt{c_{0} t^2 -1/2}}, \tau= - \frac{\pi}{2}, C=\frac{4}{{c_{0} t^2 -1/2}}, M=1, m={c_{0} t^2 +1/2}.\] Thus, (\ref{defn rho and tau}) yields $X=-\eta, Y= \xi$ and according to the results in Section \ref{sec:Constant}, it follows that
\[
e^{F}=\frac{1}{t^2\cos^2{\frac{S}{t}}},
\]
and 
$R=\frac{-2 X}{\sqrt{c_{0} t^2 -1/2}}, S= t \arctan{\sinh{u(\frac{2Y}{\sqrt{c_{0} t^2 -1/2}}}}),$ see equation (\ref{Wolf1}).

The harmonic maps between cylinders that have been constructed in \cite{Wolf2}, are similar and can be also realised as harmonic maps related to one-soliton solutions of the elliptic sinh-Gordon equation. 

Finally, there are examples in \cite{Wolf2} of harmonic maps between half-infinite cylinders. Consider the space $[0, 1] \times [1,\infty)$ equipped with the metric $\frac{dx^2 + dy^2}{y^2},$ where we identify $\{0\} \times \{y\} $ and $\{1\} \times \{y\} $ to obtain a half-infinite cylinder. Consider a harmonic map $u$ that maps the half-infinite cylinder to itself, with the boundary conditions $u(x, 1) = (x, 1)$ and $\lim_{y\rightarrow \infty} u(x, y) = \infty. $ The map $u(x, y) = x+iv(y)$ is the given harmonic, provided that $v$ satisfies $v(y)v''(y)-v'(y)^2 +1= 0$ and the initial conditions $v(1)=1, \lim_{y\rightarrow \infty}v(y)=\infty$. A non conformal solution is the one parameter family of maps \begin{equation}\label{one parameter}
v_{c}(y) =\frac{1}{\sqrt{c}} \sinh[\sqrt{c}(y-1)+\sinh^{-1}\sqrt{c}].
\end{equation} 

One can observe that the Hopf differential is given by
\begin{equation}
u_{z} \bar{u}_{z} e^{F}=\frac{1-(v')^2}{4v^2}=-\frac{c}{4}.
\end{equation}
The Beltrami coefficient is given by 
\begin{equation}
\frac{u_{\bar{z}}}{ u_{z}}=\frac{1 + v'}{1 - v'}=e^{-2 \omega}.
\end{equation}
Thus, the one parameter family of harmonic maps $u_{c}(x,y)=x+iv_{c}(y)$ can also be related to one-soliton solutions of the elliptic sinh-Gordon equation. Hence, this family can be recovered by the method presented earlier. The proof is similar, thus omitted.

\subsection{The strip model examples of Shi, Tam and Wan   \cite{S-T-W}}
This Section focuses on the explicit solution of the harmonic map problem between hyperbolic spaces given in the influential work \cite{S-T-W}, that generalizes the solutions in \cite{C-T,W}. In particular, this solution is a quasi-conformal harmonic diffeomorphism. This result can be recovered by the results in Section \ref{sec:Constant}. In fact, we shall prove that the construction given in the paper \cite{S-T-W} corresponds to the one-soliton solution of the sinh-Gordon equation.

Consider the strip model for the hyperbolic plane. In \cite{S-T-W}, the authors find a harmonic map $u=R+iS,$ satisfying the following  equations  
\begin{equation}\label{eq:Shi_eq}
 \begin{array}{c}
 R(x,y) = \alpha x + h(y) \;  \mbox{ and }  S(x,y) = g(y) \\
 a=h'(\frac{\pi}{2})\;  \mbox{  and } b=g'(\frac{\pi}{2}).
 \end{array}
 \end{equation}
   They show that $\frac{\partial R}{\partial y}= a^2 \sin^2 g $ and $\cot{g}=z,$ where 
\begin{equation}\label{shi}
\int_{0}^{z(y)}\frac{dz}{\sqrt{\alpha^2 z^4 + c^2 z^2+b^2 }}=\frac{\pi}{2}-y,
\end{equation}
and $c^2=\alpha^2+b^2+a^4.$ They extend $g,h$ to $[0,\pi]$ such that 
\[
h(y)=h(\pi)-h(\pi-y), \, g(y)=\pi-g(\pi-y),
\]
and they prove that there are suitable constants $a,b$ such that the harmonic map is a quasi-conformal harmonic diffeomorphism between the hyperbolic strips.

The same harmonic map can be recovered by the method presented in Section  \ref{sec:Constant}. More precisely, let 
\[x=X, y=Y, X_0 = 0, Y_0 = \pi/2,\]
\[
w_{1}=\frac{1}{\alpha \sqrt{2}}\sqrt{c^2 -\sqrt{c^4-4\alpha^2 b^2}}, \, w_{2}=\frac{1}{\alpha \sqrt{2}}\sqrt{c^2 +\sqrt{c^4-4\alpha^2 b^2}}.
\]
Consider 
\begin{equation}\label{eq:constants_Shi}
\begin{array}{c}
\rho=\frac{2}{\alpha \sqrt{w_2^2 -w_1^2}}, \,\tan\tau = -{\sqrt{\frac{w_2^2-1}{1-w_1^2}}}, \\
C=-\alpha^2 w_1^2,\, M=\frac{1}{w_1^2}, \,m=\frac{w_2^2}{w_1^2}, \,\Sigma=i(S-\frac{\pi}{2}), \\
\omega'_0 =0\\
\end{array}
\end{equation}
We observe that 
\begin{equation*}
\sqrt{(w_2^2-1)(1-w_1^2)}=\dfrac{a^2}{\alpha}.
\end{equation*}
Considering the choice of the parameters in \cite{S-T-W}, we have that
\[
K'(1- \dfrac{w_1^2}{w_2^2})= \alpha w_2 \dfrac{\pi}{2}=
\int_{0}^{\frac{\pi}{2}} \dfrac{d\theta}{ 1 - \dfrac{w_1^2}{w_2^2} \sin^2 \theta},
\]
where $K'$ is the imaginary  quarter period of the elliptic Jacobi functions, see equations (16.1.1) and (16.1.2) of \cite{MilThom64}. 

The relation $\omega'_0=0$ and (\ref{omega}) imply that 
\[
{cd}(v_0\vert \frac{1}{m})=   \dfrac{{cn}(v_0\vert \frac{1}{m})}{{dn}(v_0\vert \frac{1}{m})}=0\]
 Therefore, from equation (16.5.3) of the Reference  \cite{MilThom64}, we conclude that
 \[
 v_0=K( 1-\dfrac{w_1^2}{w_2^2})=\int_{0}^{\frac{\pi}{2}} \dfrac{d\theta}{ 1 -(1- \dfrac{w_1^2}{w_2^2} )  
\sin^2 \theta},
 \]
 where $K$ is the real quarter of the elliptic Jacobi functions.

From (\ref{tanhomega}), and formulas (16.8.1) and (16.20) from reference  \cite{MilThom64}, we find that
\begin{equation}\label{eq:tanh_shi}
\begin{array}{rl}
\tanh{\omega}=&
\frac{w_1}{w_2}sn(i \alpha w_2 (Y-\frac{\pi}{2})+v_0\vert 1- \frac{w_1^2}{w_2^2})=\\
=&\frac{w_1}{w_2}sn(i \alpha w_2 Y-  i K'+K\vert 1-\frac{w_1^2}{w_2^2})
=\\
=&\frac{w_1}{w_2}sn(i \alpha w_2 Y
+ i K'+K\vert 1-\frac{w_1^2}{w_2^2})
=\\
=&
dn(\alpha w_2 y\vert 1-\frac{w_1^2}{w_2^2}).
\end{array}
\end{equation} 
Similarly, we find that 
\begin{equation}\label{eq:domega_Shi}
\dfrac{\omega '(Y)}{\sqrt{C}}                                                   
=cd(i \alpha w_2 (Y-\frac{\pi}{2})+v_0\vert 1- \frac{w_1^2}{w_2^2})=
i \frac{w_{2}}{w_{1}} cs(\alpha w_2 Y \vert 1-\frac{w_1^2}{w_2^2}).
\end{equation}


Using Equation (\ref{eq:SY}),
 we find that
\[
\frac{\partial S}{\partial y}=\frac{\alpha w_2^2\, dn(\alpha w_2 y\vert  1-\frac{w_1^2}{w_2^2})}{w_2^2 + (1-w_2^2)\, sn^{2}(\alpha w_2 y\vert 1-\frac{w_1^2}{w_2^2})}.
\]
From  Equation (\ref{eq:RY}) we find that
\[
\frac{\partial R}{\partial y}=\frac{a^2 \,sn^2(\alpha w_2 y  \vert  1-\frac{w_1^2}{w_2^2})}{w_2^2 + (1-w_2^2)\, sn^{2}(\alpha w_2 y \vert  1-\frac{w_1^2}{w_2^2})}.
\]

Now we shall compare the above results to the results given in Reference \cite{S-T-W}.
 
 Note that equation (\ref{shi}) can be written as 
\[
\int_{0}^{z(y)}\frac{dz}{\alpha\sqrt{(z^2 + w_1^2)(z^2 + w_2^2) }}=\frac{\pi}{2}-y.
\]
This equation can be solved by using elliptic integrals, see for example \cite[p. 217]{Arm-Eber}. The solution is 
\[
z=w_1 sc\left(\alpha w_2 (\frac{\pi}{2}-y \vert  1-\frac{w_1^2}{w_2^2}\right).
\]
Taking into account that  $K'= \alpha w_2 \dfrac{\pi}{2}$ we find that
\[
z=w_2 cs(\alpha w_2 y \vert  1-\frac{w_1^2}{w_2^2}).
\]
Thus,
\[
S= \cot^{-1}{z}=\cot^{-1}{\left( w_2 cs(\alpha w_2 y \vert  1-\frac{w_1^2}{w_2^2}) \right)}.
\]
It follows that 
\begin{equation}
\frac{\partial R}{\partial y}= h'(y) = a^2 \sin^2{S}=\frac{a^2 \,sn^2(\alpha w_2 y  \vert  1-\frac{w_1^2}{w_2^2})}{w_2^2 + (1-w_2^2)\, sn^{2}(\alpha w_2 y \vert  1-\frac{w_1^2}{w_2^2})}.
\end{equation}
Similarly,
\begin{equation}
\begin{split}
\frac{\partial S}{\partial y}&= -\frac{1}{1+w_2^2 cs^2(\alpha w_2 y \vert  1-\frac{w_1^2}{w_2^2})} w_2^2 \alpha  (cs)'(\alpha w_2 y \vert  1-\frac{w_1^2}{w_2^2})\\
& = \frac{\alpha w_2^2\, dn(\alpha w_2 y\vert  1-\frac{w_1^2}{w_2^2})}{w_2^2 + (1-w_2^2)\, sn^{2}(\alpha w_2 y\vert 1-\frac{w_1^2}{w_2^2})}.
\end{split}
\end{equation}

This result is identical with the result above.
Then, (\ref{eq:relation_eFRY}) and (\ref{eq:rel}) coincide with the corresponding equations (4.3) in \cite{S-T-W}. More precisely, we find that
\begin{equation*}
 \dfrac{\partial R}{\partial Y}=-\dfrac{4  \sin\, \tau\, \cos\, \tau}{\alpha \rho^2}  \dfrac{1}{\Phi}=a^2 \sin^2 S
\end{equation*}
and
\begin{equation*}
\begin{array}{rl}
\left(\dfrac{\partial S}{\partial Y}\right)^2=&
-\dfrac{16 {\tan^2 \,\tau }}{\left(\tan^2 \, \tau +1\right)^2
   \alpha ^2 \rho ^4}\dfrac{1}{\Phi^2}
   +\dfrac{4 \left({\tan^2\,\tau
   }-1\right)}{\left({\tan\,\tau  }^2+1\right) \rho ^2 \Phi
   }+\alpha ^2\\
   =&\alpha ^2 + (b^2+a^4-\alpha^2)\sin^2{S}-a^4 \sin^4{S}.
   \end{array},
\end{equation*}
Also,
\begin{equation}\label{eq:S}
\begin{array}{rl}
S-\frac{\pi}{2}=&
i
{
\mathrm{arctanh} \left( w_1 \dfrac{\omega'(Y)}{\sqrt{C} }\right),
}
\\
\end{array}
\end{equation}
\begin{equation*}
e^F= \frac{1}{\sin^2{S}}, 
\end{equation*}
and using equation (\ref{eq:domega_Shi}) we find
\begin{equation*}
S=\cot^{-1}\left( w_2 cs(\alpha w_2 Y\vert 1-\frac{w_1^2}{w_2^2})  \right),
\end{equation*}
and this result coincides with the result in \cite{S-T-W}.

\subsection{The upper half-space example of Li and Tam \cite{Li-Ta2}}\label{L-T}
This Section focuses on the explicit solution of the harmonic map problem in the influential work \cite{Li-Ta2}. This result can also be recovered by the results in Section \ref{sec:Constant}, as shown below. 

Consider the solution 
\[
\omega(\zeta,\bar{\zeta})=-\log{\tanh{\xi}}
\]
of the equation
\[
\Delta \omega = 2 \sinh{2\omega}.
\]
Let $a>0$. Then, the equation
\[
\frac{u_{\bar{\zeta}}}{u_{{\zeta}}}=e^{-2\omega}
\]
admits the solutions 
\[
u(\zeta,\bar{\zeta})=(\frac{2\eta}{a}, -2\sinh \frac{2\xi}{a}).
\]
If $z=-\frac{2i}{a}\zeta,$ then 
\[
u(z,\bar{z})=(x, \frac{1}{a}\sinh{ay}).
\]
The hyperbolic metric that corresponds to the target $N$ is $e^{F(R,S)}=\frac{1}{S^2},$ see Section \ref{sec:Constant} for more details. This is a family of harmonic maps between hyperbolic spaces, that has been studied in \cite{Li-Ta2}. 

Note that this result can also be obtained by the general result in Section \ref{sec:Constant}, by considering the initial conditions $u(x_0,y_0)=(x_0,\frac{1}{a}\sinh{ay_{0}}),$ $ \omega_{0}=-\log\tanh{y_0}, \omega_0^{'} = \frac{1}{\cosh{2y_0}}$. In this case $K_{N}=-1, \rho=1, C= 0, \tau= -\frac{\pi}{2}$ and one can observe that the limit case of the solution in Section \ref{sec:Constant}, taking $C\rightarrow 0, \tau \rightarrow  -\frac{\pi}{2},$ converges pointwise to the explicit solution obtained in \cite{Li-Ta2}.

\section{B\"acklund Transform of the sinh-Gordon equation}\label{sec:Backlund}

In this section we discuss a B\"acklund transform of the sinh-Gordon equation, by applying the methods discussed in the previous sections. A new harmonic map is provided by using the B\"acklund transform of a solution to the sine-Gordon equation. 

The \textit{B\"acklund transformation} is a system of first order partial differential equations relating the solution of a PDE   (in our case the sinh-Gordon PDE ) to the solution of  another PDE (in our case the  PDE (\ref{eq:theta_zeta_barzeta})). Then, the one solution is said to be the \textit{B\"acklund transform} of the other. 
 
In what follows we shall prove the following proposition. 
\begin{proposition}
The  equation  
\begin{equation}\label{eq:theta_zeta_barzeta}
      2 i \theta_{ \zeta \bar{\zeta}}  =
        e^{-F}
        \left(     \left( F_u^2-F_{uu}\right) e^{2 i \theta}   -  \left( F_{\bar{u}}^2-F_{\bar{u}\bar{u}}\right) e^{-2 i \theta} \right)
\end{equation}
is the B\"acklund transform of the sinh-Gordon equation
\begin{equation}\label{eq:sinh-G_zeta}
\omega_{\zeta \bar{\zeta}}=-\frac{K_N}{2} \sinh 2 \omega, 
\end{equation}
where \[ K_N=K_N(u,\bar{u})
=-\dfrac{1}{2} \mathrm{e}^{-F} \Delta F\] is the Gauss curvature of the metric
 \[
h= \mathrm{e}^{F(u,\bar{u})}   |du|^2.
 \]
 In addition, the system of first order partial differential equations that relates the two solutions is (\ref{eq:Back_system}).
 \end{proposition}

Using the specific coordinate system (\ref{eq:specific}), equations (\ref{eq:exp_v}) and (\ref{eq:exp_vbar}) can be written as follows\,:
  \begin{equation}\label{eq:exp_v_spec}
  e^{v- \frac{F}{2}}\,= {u}_\zeta
  , \quad
   e^{-v - \frac{F}{2}}\,= \bar{u}_\zeta, \quad
 e^{\bar{v}- \frac{F}{2}}\, ={\bar{u}}_{\bar{\zeta}}
 , \quad
   e^{-\bar{v} - \frac{F}{2}}\,= {u}_{\bar{\zeta}}, 
   \quad
      v=\omega+i \theta.
 \end{equation} 
Consider the relation
\[
\left(u_{\bar{\zeta}} \right)_{\zeta}=\left(u_{\zeta} \right)_{\bar{\zeta}}\,.
\]
Then $v$ satisfies the linear partial differential equation
\begin{equation}\label{eq:comp}
\left( e^{v- \frac{F}{2}}\right)_{\bar{\zeta}}=
\left(e^{-\bar{v} - \frac{F}{2}}\right)_{\zeta}.
\end{equation}
Equating real and imaginary parts, we observe that (\ref{eq:comp}) can be written as a system of first order nonlinear differential equations:
\begin{equation}\label{eq:Back_system}
\begin{array}{rl}
\omega_{\xi}-\theta_{\eta}=& \dfrac{1}{2} \tanh\, \omega\, \dfrac{\partial F}{\partial \xi}
\\   \\
\omega_{\eta}+\theta_{\xi}=& \dfrac{1}{2} \coth\, \omega\, \dfrac{\partial F}{\partial \eta}
\end{array} .
\end{equation}
By the chain rule we have
\begin{equation*}
\dfrac{\partial F}{\partial \xi}=F_u\, u_{\xi} + F_{\bar{u}}\, {\bar{u}}_{\xi} , \; 
\dfrac{\partial F}{\partial \eta}=F_u\, u_{\eta} + F_{\bar{u}}\, {\bar{u}}_{\eta}.
\end{equation*}
In these equations, we can replace the partial derivatives of the function $u$ by the relations given in equation (\ref{eq:exp_v_spec}). One can check that the function $\theta(\xi,\eta)$ satisfies equation  (\ref{eq:theta_zeta_barzeta}) and $\omega(\xi,\eta)$ satisfies the sinh-Gordon type equation (\ref{eq:sinh-G_zeta}). Therefore the function $\theta$ is indeed the B\"acklund transform of the function $\omega$. Let us outline the calculations. 

From equation (\ref{eq:def_w}) and taking into account (\ref{eq:harmonic_map}), it follows that
\begin{equation}\label{backlund proof}
2v_{\bar{\zeta}}=F_{\bar{u}}\bar{u}_{\bar{\zeta}}-F_{u}u_{\bar{\zeta}},
\end{equation}
and thus, using the chain rule, it follows that
\begin{equation}
    \begin{split}
        2 v_{ \zeta \bar{\zeta}} =& 2 \omega_  { \zeta \bar{\zeta}} +2 i \theta_{ \zeta \bar{\zeta}}  =\\
        =&F_{ u \bar{u}} \left( |u_\zeta|^2 - | u_{\bar{\zeta}}|^2\right) +  \left( F_u^2-F_{uu}\right) u_z u_{\bar{\zeta}} -  \left( F_{\bar{u}}^2-F_{\bar{u} \bar{u}}\right) {\bar{u}}_{\zeta} {\bar{u}}_{\bar{\zeta}} .
        \end{split}
\end{equation}
Therefore,
\[
2\omega_{\zeta\bar{\zeta}} = F_{u\bar{u}} \left(|u_{\zeta}|^2- |u_{\bar \zeta}|^2 \right) ,
\]
and
\[
 2 i \theta_{ \zeta \bar{\zeta}}  =
          \left( F_u^2-F_{uu}\right) u_{\zeta} u_{\bar{\zeta}} -  \left( F_{\bar{u}}^2-F_{\bar{u} \bar{u}}\right) {\bar{u}}_{\zeta} {\bar{u}}_{\bar{\zeta}} .
\]
Taking into consideration (\ref{eq:exp_v}) and  (\ref{eq:exp_vbar}), one can calculate the following relations

\begin{equation*}\label{eq:exptheta}
e^{2 i \theta} = e^{v-\bar{v}}= e^{F}  u_{\zeta}  u_{\bar{\zeta}},  \quad 
e^{-2 i  \theta} = e^{-(v-\bar{v})}= e^{F}  \bar{u}_{\zeta}  \bar{u}_{\bar{\zeta}},
\end{equation*}   
and
\begin{equation*}\label{eq:sinhtheta}
\sin(2 \theta)=  \dfrac{ e^{F}}{2 i }  \left(u_{\zeta}  u_{\bar{\zeta}}
 -\bar{u}_{\zeta}  \bar{u}_{\bar{\zeta}}\right),   \;
\cos(2 \theta)=  \dfrac{ e^{F}}{2}  \left(
u_{\zeta}  u_{\bar{\zeta}}
+\bar{u}_{\zeta}  \bar{u}_{\bar{\zeta}}\right).
\end{equation*}
From the above it follows that 
 \begin{equation*}\label{eq:theta_zbarz}
    \begin{split}
      2 i \theta_{ \zeta \bar{\zeta}}  =e^{-F}
        \left(     \left( F_u^2-F_{uu}\right) e^{2 i \theta}   -  \left( F_{\bar{u}}^2-F_{\bar{u}\bar{u}}\right) e^{-2 i \theta} \right).
        \end{split}
\end{equation*}

Similarly, as already mentioned, the following holds true:
\[
2\omega_{\zeta\bar{\zeta}} = F_{u\bar{u}} \left(|u_{\zeta}|^2- |u_{\bar \zeta}|^2 \right) 
=2 F_{u \bar{u}} e^{-F}\sinh 2 \omega= - K_N \sinh 2 \omega,
\]
thus,
 \begin{equation*}\label{eq:sinh-Gordon_specific}
    \omega_{\zeta\bar{\zeta}}= - \frac{K_N}{2}    \sinh 2 \omega.
 \end{equation*}    

One interesting remark arising directly from the above, is the following.
\begin{remark}
Consider the case of the hyperbolic upper half-plane metric, where
\[
F=2 \log \left( \dfrac{2i}{u - \bar{u}} \right).
\]
Then, $F_{u}=ie^{\frac{F}{2}}, F_{\bar{u}}=-ie^{\frac{F}{2}}$ and thus $F_{\zeta}=-2\sin{\theta}\cosh{\omega}+2i\cos{\theta}\sinh{\omega}.$

Then, the function $\omega$ satisfies the sinh-Gordon equation
\[
\omega_{\zeta \bar{\zeta}}=\frac{1}{2} \sinh(2 \omega)
\]
and the function $\theta$ satisfies the sine-Gordon equation
 \[
\theta_{\zeta \bar{\zeta}}=-\frac{1}{2} \sin(2 \theta).
\]
From (\ref{eq:Back_system}) we find that these two functions are related by the transform
\begin{equation}\label{eq:Backlund_upper}
\begin{array}{rl}
\omega_\xi- \theta_\eta&= -2 \sinh\,\omega \sin\, \theta\\
\omega_\eta+\theta_\xi&=  -2 \cosh\,\omega \cos\, \theta
\end{array}.
\end{equation}
This is a B{\"a}cklund transform and it provides a connection between the solutions of an elliptic sinh-Gordon equation and an elliptic sine-Gordon equation. Thus, one can obtain solutions of the elliptic sinh-Gordon equation by the known solutions of the sine-Gordon equation and vice versa.
\end{remark}

Taking into account the last Remark, one can construct a solution of the elliptic sinh-Gordon equation by a one-soliton solution of the sine-Gordon equation.  Consider for example the solution \[
\theta(\xi,\eta)=\arcsin{\tanh{2\xi}}
\] 
of the equation 
 \begin{equation}\label{theta}
 \theta_{\zeta \bar{\zeta}}=-\frac{1}{2} \sin(2 \theta).
 \end{equation}
Then, substituting (\ref{theta}) into (\ref{eq:Backlund_upper}), it follows that 
\[
\omega(\xi,\eta)=2 \text{arctanh}{\frac{2\eta}{\cosh{(2\xi)}}}, \,\, \eta> \frac{1}{2}\cosh{2\xi},
\]
is the B{\"a}cklund transform  of $\theta$ and satisfies the elliptic sinh-Gordon equation
\[
\omega_{\zeta \bar{\zeta}}=\frac{1}{2} \sinh(2 \omega),
\]
i.e. $\omega$ is not a one-soliton solution of the sinh-Gordon equation, in coordinates $(\xi, \eta)$. 

A harmonic map $u=R+iS$ associated to the solution $\omega$ is 
\[
u(\xi,\eta)=\left(\eta^2\tanh{2\xi}+\frac{\xi}{2}\right)+ i\left(\frac{\eta^2}{\cosh{2\xi}}-\frac{\cosh{2\xi}}{4}\right), \, \eta> \frac{1}{2}\cosh{2\xi}.
\]
Finally, observe that the metric on the target is given by
\[
\begin{array}{rl}
e^{F(R,S)}= &\dfrac{1}{u_{\zeta} \bar{u}_{{\zeta}}}=\frac{4}{\left(  \frac{\eta^2}{\cosh{2\xi}}-\frac{\cosh{2\xi}}{4}\right)^2}
= \frac{1}{S^2},
\end{array}
\]
i.e. the hyperbolic metric of the upper-half plane, as it was expected.

\section{Perspectives for future investigation}\label{sec:prespectives}

An  application of the methods introduced in this paper is to start from the known solutions of the sinh-Gordon equation and consequently generate new solutions of the harmonic diffeomorphisms to a constant curvature manifold. This project is currently under investigation.

 The classification of the solutions of the sinh-Gordon equation is still a difficult  open problem. The solution of this problem would lead to the classification of the harmonic diffeomorphisms to a constant curvature surfaces, and by conformally equivalent mappings to the classification of all harmonic diffeomorphisms.

It would be interesting to generalize the above results in higher dimensions. Furthermore, theoretical results could be obtained by applying the theory of Beltrami equation in order to study harmonic maps between surfaces. 

On the other hand, one can apply maximum principle arguments and deduce theoretical results on harmonic diffeomorphisms. There is also a connection with the theory of constant mean curavature surfaces, implied by the sinh-Gordon equation. The case of the flat complex plane is still of great interest. More generally, it is interesting to find solutions (for example group invariant solutions) of the sinh-Gordon equation and then obtain families of harmonic maps. On the other hand, the known formulas for harmonic maps can provide examples of solutions to the elliptic sinh-Gordon and elliptic sine-Gordon equation. Note that the B{\"a}cklund transform obtained relates known solutions of the first equation to solutions of the second one.

\section*{Acknowledgement}The authors would like to thank Professor Michael Wolf for his valuable suggestions to take into consideration, which relate our results to the work of many authors. These remarks led us to radically modify  the initial version of this article. The authors would also like to thank the anonymous referee for the valuable comments and suggestions, which led us to clarify and simplify many details of the initial submitted version of the article.

\appendix
\section{Elliptic Functions summary}\label{app:elliptic}
In this Appendix, the definitions and properties of the elliptic integrals and Jacobi functions, used in this paper, are presented for clarity reasons. 
The formulation which are used in this paper,  is taken from \cite{MilThom64}.
The elliptic integral of the first kind $F(\phi \vert n)$  and the Jacobi elliptic function $sn(v\vert n)$  are defined by the formula
\begin{equation}\label{eq:Jacobi_ampli}
F(\phi\vert n) = v=
\int_0^x \, \dfrac{dt}{   \sqrt{ \left( 1-t^2\right)\left( 1 -n t^2\right)}}=
sn^{-1}(x\vert n),
\end{equation}
where
\[
 sn(v \vert n)=x=\sin \phi , \quad cn(v \vert n)=\cos \phi ,
  \quad dn(v \vert n)= \sqrt{1-n  \sin^2\, \phi}
\]
%
%

The Jacobi functions satisfy the following well known relations
\begin{equation}\label{eq:relations_Jacobi}
1- n \,sn^2(v \vert n) = dn^2(v \vert n)
\; \mbox{and}\;  sn^2(v\vert n)+cn^2(v \vert n)=1,
\end{equation}
If $p,q,r$ are any of the letters $s,d,c,n$  then  the other Jacobi functions are defined by: 
\begin{equation}\label{eq:All_Jacobi}
pq(v \vert n)=\frac{pr(v \vert n)}{qr(v \vert n)},\quad pp(v \vert n)=1.
\end{equation}

Equation (\ref{eq:mod_Sinh_G}) can be written as the sinh-Gordon equation
\begin{equation}\label{eq:sinh_G_app}
 \dfrac{\omega'(Y)}{\sqrt{ 1- (m-1) \sinh^2 \omega} } =\epsilon \sqrt{C}, \; \epsilon=\pm 1. 
\end{equation}
If we put 
\[
v(Y)= \sqrt{m} \tanh \omega(Y),
\] 
then (\ref{eq:sinh_G_app}) can be rewritten as
\[
\dfrac{v'(Y)}{ \sqrt{ (1- v^2) ( 1- \dfrac{1}{m} v^2)}}= \epsilon \sqrt{Cm} , \quad \epsilon=\pm 1
\]
Therefore, after integration, equation (\ref{eq:sinh_G_app}) yields
\begin{equation}\label{A3}
\int^{\sqrt{m}\tanh{\omega}}_{  \sqrt{m}\tanh{\omega_0}}
\dfrac{dv}{ \sqrt{ (1- v^2) ( 1- \dfrac{1}{m} v^2)}}=\epsilon \sqrt{Cm}\left(Y-Y_0\right).
\end{equation}
The left hand side of the equation (\ref{A3}) can be replaced by definition of the Jacobi elliptic functions (\ref{eq:Jacobi_ampli}), and thus $\omega$ can be found explicitly.  More precisely, we have that 
\[
v-v_0=sn^{-1}( \sqrt{m} \tanh\omega \vert \dfrac{1}{m})- 
sn^{-1}( \sqrt{m} \tanh\omega_0 \vert \dfrac{1}{m})=\epsilon \sqrt{Cm} (Y-Y_0)
\]
or
\begin{equation}\label{eq:tanhw}
\tanh{\omega}=\dfrac{1}{\sqrt{m}} sn({\epsilon} \sqrt{Cm}(Y-Y_0)+v_0 \vert \frac{1}{m}),
\end{equation}
where 
\begin{equation}
v_0 =sn^{-1}(\sqrt{m}\tanh{\omega_0}\vert \frac{1}{m})= sd^{-1}(\sqrt{m}\sinh{\omega_0}\vert \frac{1}{m}).
\end{equation}
From the relations between the Jacobi elliptic functions (\ref{eq:relations_Jacobi})
 it follows that
\begin{equation}\label{eq:sinhw}
\sinh{\omega}=
\dfrac{1}{\sqrt{m}} sd({\epsilon}\sqrt{Cm}(Y-Y_0)+v_0 \vert \frac{1}{m}),
\end{equation}
and hence
\begin{equation}
\label{eq:coshw}
\cosh{\omega}= nd(\sqrt{Cm}(Y-Y_0)+v_0 \vert \frac{1}{m}).
\end{equation}
Thus
\begin{equation}
{\omega}=\log{\frac{sd(\epsilon \sqrt{Cm}(Y-Y_0)+v_0 \vert \frac{1}{m})+\sqrt{m}\, nd(\epsilon \sqrt{Cm}(Y-Y_0)+v_0 \vert \frac{1}{m})}{\sqrt{m}} }.
\end{equation}
Furthermore, 
\begin{equation}\label{eq:omega_pr}
\begin{array}{rl}
\dfrac{\omega '(Y)}{\sqrt{C}}                                                   
=&\epsilon  cd(\epsilon \sqrt{Cm}(Y-Y_0)+v_0 \vert \frac{1}{m})\\
=&  \epsilon  sn(\epsilon \sqrt{Cm}(Y-Y_0)+v_0 +K(\dfrac{1}{m})  \vert \frac{1}{m}).
\end{array}
\end{equation}

From equation (\ref{eq:SY}), after some elementary calculations, it follows that
\begin{equation}\label{eq:SY_app}
S_Y= -\alpha \dfrac{\cosh \omega  \sinh \omega}{\cos^2 \tau  + \sinh^2 \omega}= \dfrac{\alpha}{C M}
\dfrac{\omega''}{1-  \left(\dfrac{\omega'(Y)}{\sqrt{CM}}\right)^2 }.
\end{equation}
Using the formulas (\ref{eq:coshw}) and (\ref{eq:sinhw})
one can find that
\begin{equation}\label{eq:SY-tel}
S_Y=
-\dfrac{\alpha}{\sqrt{m}}   
\dfrac{  nd(v\vert \frac{1}{m})  sd(v\vert \frac{1}{m}) }{  \cos^2 \,\tau  +  \dfrac{1}{m}  nd^2(v \vert \frac{1}{m})},
\end{equation}
where 
\[
v=\epsilon \sqrt{Cm}(Y-Y_0)+v_0,
\]
and  we can verify the implicit formula 
(\ref{eq:S}).
Formula (\ref{eq:S_fin}) can be obtained by integrating the right hand side of equation (\ref{eq:SY_app}).
 
From equation (\ref{eq:RY}), by similar algebraic calculations we find
\[
R_{Y}=-\dfrac{ \alpha \tan \tau  } { \cosh^2 \omega ( 1+ \tan^2 \tau \tanh^2\omega )}=
-\dfrac{ \alpha \tan \tau \, dn^2(v\vert \frac{1}{m}) } {  ( 1+ \dfrac{1}{m}\tan^2 \tau \, sn^2(v\vert \frac{1}{m}) )}.
\]
We also find
 \begin{equation}\label{eq:RY_app}
\begin{array}{rl}
R_Y=& -	\dfrac{(m-1) \,\alpha \, \sin\, \tau\,  \cos\,\tau  }{M}
\dfrac{1}{ 1-\dfrac{1}{M}
\left(  \dfrac{\omega'(Y)}{\sqrt{C}}\right)^2}\\
=&
 -	\dfrac{(m-1) \,\alpha \, \sin\, \tau\,  \cos\,\tau  }{M}
\dfrac{1}{ 1-\dfrac{1}{M}
  sn^2(\epsilon \sqrt{Cm}(Y-Y_0)+v_0 +K(\dfrac{1}{m})  \vert \frac{1}{m})}.
\end{array}
\end{equation}

By definition  of the third kind  elliptic integral, we have
\[
\Pi\left( n; \,u \,\vert\, \ell  \right)\, =\,
\int_{0}^{u}\,
\dfrac{dw}{  \sqrt{  ( 1- n\;\rm{sn}^2 (w\vert \ell))}}.
\]
Integrating (\ref{eq:RY_app}), we find
\[
\begin{array}{l}
R-R_0=
\alpha (X-X_0) -  \\
- \epsilon \sqrt{ C m} 	\dfrac{(m-1) \,\alpha \, \sin\, \tau\,  \cos\,\tau  }{M}
\left(
\Pi\left( \dfrac{1}{M}; \,v +K \,\vert\,\dfrac{1}{m}  \right)-
\Pi\left(
 \dfrac{1}{M}; \,v_0 +K \,\vert\,\dfrac{1}{m} 
\right)
\right),
\end{array}
\]
thus we can obtain formula (\ref{eq:R_Pi}).
 


\begin{thebibliography}{99}
\bibitem{A-R}A. Antonio and S. Rabah:
\emph{Harmonic diffeomorphisms between domains in the Euclidean 2-sphere}.
Comment. Math. Helv. 89 (2014), no. 1, 255-271\,.
\bibitem{C-T}H-I Choi and A. Treibergs:
\emph{New examples of harmonic diffeomorphisms of the hyperbolic plane onto itself}. Manuscripta Math. 62 (1988), no. 2, 249–256\,.

\bibitem{DW}D. Dumas and M. Wolf:
\emph{Polynomial cubic and differentials and convex polygons in the projective plane}. Geom. Funct. Anal. 25 (2015), no. 6, 1734–1798\,.

\bibitem{Arm-Eber} J. V. Armitage and W. F. Eberlein: 
\emph{Elliptic Functions}. London Mathematical Society, Student Texts 67, Cambridge University Press (2006).


\bibitem{Joa}P. Joaquín:
\emph{Sinh-Gordon type equations for CMC surfaces}. 145-156, Ed. Univ. Granada, Granada, 2011\,. 
\bibitem{J}J. Jost:
\emph{Riemannian Geometry and Geometric Analysis}. (Third Edition)
Springer (2002).
\bibitem{F-P} A. S. Fokas and B. Pelloni:
\emph{The Dirichlet-to-Neumann map for the elliptic sine-Gordon equation}. Nonlinearity 25 (2012), no. 4, 1011-1031\,. 

\bibitem{G-R-S-Y}V. Gutlyanskii, V. Ryazanov, U. Srebro, and E. Yakubov:
\emph{The Beltrami Equation. A geometric approach}. Developments in Mathematics, 26. Springer, New York, 2012. xiv+301 pp\,.

\bibitem{Han}Z-C. Han:
\emph{Remarks on the geometric behavior of harmonic maps between surfaces}. Elliptic and parabolic methods in geometry (Minneapolis, MN, 1994), 57–66, A K Peters, Wellesley, MA, 1996\,.

\bibitem{H-T-T-W}Z-C. Han, L-F. Tam, A. Treibergs and T. Wan:
\emph{Harmonic maps from the complex plane into surfaces with nonpositive curvature}.
Comm. Anal. Geom. 3 (1995), no. 1-2, 85–114\,.


\bibitem{Hitchin}N. J. Hitchin:
\emph{The self-duality equations on a Riemann surface}. 
Proc. London Math. Soc. (3) 55 (1987), no. 1, 59–126\,. 
\bibitem{H}G. Hwang:
\emph{The elliptic sinh-Gordon equation in the half plane}. 
J. Nonlinear Sci. Appl. 8 (2015), no. 2, 163-173\,.
\bibitem{Kal}D. Kalaj:
\emph{Rad\'{o}-Kneser-Choquet theorem for harmonic mappings between surfaces}. 
Calc. Var. Partial Differential Equations 56 (2017), no. 1, Art. 4, 12 pp\,.
\bibitem{K}K. Kenmotsu:
\emph{Surfaces with constant mean curvature}. Translations of Mathematical Monographs, 221. American Mathematical Society, 2003\,. 
\bibitem{Li-Ta2}P. Li and L-F. Tam:
\emph{Uniqueness and regularity of proper harmonic maps}. Ann.
of Math. (2) 137 (1993), 167-201\,.

\bibitem{Loftin}J. C. Loftin: 
\emph{Affine spheres and convex RPn-manifolds}. Amer. J. Math. 123 (2001), no. 2, 255–274\,.


\bibitem{M2}V. Markovic:
\emph{Harmonic diffeomorphisms and conformal distortion of Riemann surfaces}. Comm. Anal. Geom. 10 (2002), no. 4, 847-876\,.
\bibitem{M}V. Markovic:
\emph{Harmonic maps and the Schoen conjecture}. J. Amer. Math. Soc. 30 (2017), no. 3, 799-817\,.
\bibitem{MilThom64}
L. M. Milne-Thomson:\emph{ Jacobian Elliptic and Theta Functions} Chapter 16 and
\emph{Elliptic Integrals} Chapter 17 in \emph{Handbook of Mathematical Functions}, Editors M. Abramowitz and I. Stegun, 
\emph{Applied Mathematics Series 55, National Bureau Standards}, 10th ed., 1964
\bibitem{Minsky}Y. N. Minsky:
\emph{Harmonic maps, length and energy in Teichm\"{u}ller space}. J. Differential Geom. 35 (1992), no. 1, 151–217\,. 
\bibitem{P-H}C. Pascal and R. Harold:\emph{
Construction of harmonic diffeomorphisms and minimal graphs}.
Ann. of Math. (2) 172 (2010), no. 3, 1879-1906\,. 
\bibitem{P-R}Y. Pinchover, R. and J. Rubinstein:
\emph{An introduction to partial differential equations}. Cambridge University Press, Cambridge, 2005\,.
\bibitem{S-Y}R. Schoen and S-T. Yau:
\emph{On univalent harmonic maps between surfaces}. 
Invent. Math. 44 (1978), no. 3, 265-278\,. 
\bibitem{S-T-W}Y. Shi, L-F. Tam and T. Y. H. Wan:
\emph{Harmonic Maps on Hyperbolic Spaces with Singular Boundary
Value}. Journal of Differential Geometry 51 (1999), no. 3, 551-600\,.
\bibitem{W} J. Wang: 
\emph{The heat flow and harmonic maps between complete manifolds}. J. Geom. Anal. 8 (1998), no. 3, 485-514\,.

\bibitem{Wente} H. C. Wente: 
\emph{Counterexample to a conjecture of H. Hopf}. Math. 121 (1986), no. 1, 193-243\,.

\bibitem{Wolf1}M. Wolf:
\emph{The Teichm\"{u}ller theory of harmonic maps}. J. Differential Geom. 29 (1989), no. 2, 449–479\,. 
\bibitem{Wolf2}M. Wolf:
\emph{Infinite energy harmonic maps and degeneration of hyperbolic surfaces in moduli space}. J. Differential Geom. 33 (1991), no. 2, 487–539\,. 
\bibitem{Wolf3}M. Wolf:
\emph{High energy degeneration of harmonic maps and rays in Teichm\"{u}ller space}. Topology 30 (1991), no. 4, 517–540\,. 
\bibitem{X-A}C. Xingdi and F. Ainong:
\emph{A note on harmonic quasiconformal mappings}. J. Math. Anal. Appl. 348 (2008), no. 2, 607-613\,. 
\end{thebibliography}
\end{document}